\RequirePackage{ifpdf}
\ifpdf 
\documentclass[pdftex]{sigma}
\else
\documentclass{sigma}
\fi

\numberwithin{equation}{section}

\newcommand{\e}[1]{{(\ref{#1})}}
\newcommand{\eq}[1]{{equation~(\ref{#1})}}
\newcommand{\es}[2]{{(\ref{#1}) and (\ref{#2})}}
\newcommand{\eqs}[2]{{equations (\ref{#1}) and (\ref{#2})}}

\newcommand{\mb}[1]{{\mbox{${#1}$}}}

\newcommand{\RHS}{{\rm RHS}}

\newcommand{\eps}{\varepsilon^{}}
\newcommand{\ocomma}{\stackrel{\circ}{,}}
\newcommand{\starcomma}{\stackrel{*}{,}}
\newcommand{\veee}{\vee}
\newcommand{\bigveee}{\bigvee{}}
\newcommand{\cA}{{\cal A}}
\newcommand{\cR}{{\cal R}}
\newcommand{\casimir}{C}
\newcommand{\mult}{m}
\newcommand{\mm}{m}
\newcommand{\nn}{n}
\newcommand{\omegaone}{\varpi}

\newcommand{\odd}{{\rm odd}}
\newcommand{\even}{{\rm even}}
\newcommand{\cl}{{\rm cl}}
\renewcommand{\hbar}{\mathchar'26\mkern-9mu h}
\newcommand{\Hf}{\frac{1}{2}}
\newcommand{\papar}[1]{
\frac{\stackrel{\raise.2ex\hbox{$\leftarrow$}}{\partial^{r}}   }
{   \partial {#1}}  }
\newcommand{\papal}[1]{
\frac{\stackrel{\lower.3ex \hbox{$\rightarrow$}}{\partial^{\ell}}   }
{   \partial {#1}}  }
\newcommand{\larrow}[1]{\stackrel{\rightarrow}{#1}}

\begin{document}

\allowdisplaybreaks

\renewcommand{\thefootnote}{$\star$}

\renewcommand{\PaperNumber}{036}

\FirstPageHeading

\ShortArticleName{Three Natural Generalizations of Fedosov Quantization}

\ArticleName{Three Natural Generalizations\\ of Fedosov Quantization\footnote{This paper
is a contribution to the Special Issue on Deformation
Quantization. The full collection is available at
\href{http://www.emis.de/journals/SIGMA/Deformation_Quantization.html}{http://www.emis.de/journals/SIGMA/Deformation\_{}Quantization.html}}}

\Author{Klaus BERING}

\AuthorNameForHeading{K.~Bering}

\Address{Institute for Theoretical Physics {\rm \&} Astrophysics,
Masaryk University,\\ Kotl\'a\v{r}sk\'a 2, CZ-611 37 Brno, Czech Republic}

\Email{\href{mailto:bering@physics.muni.cz}{bering@physics.muni.cz}}

\ArticleDates{Received May 19, 2008, in f\/inal form February 14,
2009; Published online March 25, 2009}

\Abstract{Fedosov's simple geometrical construction for deformation quantization of
symplectic manifolds is generalized in three ways without introducing new
variables: (1)~The base manifold is allowed to be a supermanifold. (2)~The
star product does not have to be of Weyl/symmetric or Wick/normal type. (3)~The initial geometric structures are allowed to depend on Planck's constant.}


\Classification{53D05; 53D55; 58A15; 58A50; 58C50; 58Z05}

\section{Introduction}
\label{secintro}

The aim of this paper is to give natural generalizations of Fedosov's {\em
simple geometrical construction for deformation quantization} on a symplectic
manifold \mb{M} \cite{f85,f86,f94,f96}. In short, the term ``deformation
quantization'' refers to the construction of an associative \mb{*} product that
is an expansion in some formal parameter \mb{\hbar}, and whose leading behavior
is controlled by a~geo\-metric structure \mb{\mult^{ij}}, which usually (but not
always) is a symplectic structure \cite{bffls78}. {}Fedosov quantization, in
its most basic form, is a deformation quantization recipe that relies on yet
another geometric input in form of a compatible torsion-free tangent bundle
connection \mb{\nabla}. It is an important feature of Fedosov's \mb{*}
product that it is a~{\em global} construction, which applies to a manifold as
a~whole, and not just a local neighborhood or vector space. Another
characteristic feature that sets Fedosov's approach apart from most other
methods is the assemblage of a globally well-def\/ined, f\/lat connection \mb{D}.
We discuss in this paper the following three natural generalizations of
Fedosov's original construction:

\begin{enumerate}\itemsep=0pt
\item
We let the base manifold be an arbitrary smooth supermanifold. Previous works
on supe\-ri\-zing Fedosov's construction use Batchelor's theorem (which in turn
relies on the existence of a partition of unity) to give a non-canonical
identif\/ication of the supermanifold with the ringed space of sections of the
exterior algebra of a vector bundle, so that the bosonic and fermionic
variables are associated with base and f\/iber directions in this vector bundle,
respectively, see e.g., \cite{bordemann96,rothstein90}. On the contrary, we
will treat bosons and fermions on equal footing, and will locally allow the use
of arbitrary coordinate systems. This of course requires that we pay special
attention to sign factors, but it will be more elementary~-- and more
general~-- in the sense that we will not need any major results from the
theory of supermanifolds.
\item
We do not assume that \mb{\mult^{ij}} is antisymmetric. This is sometimes
referred to as the \mb{*} product does not have to be of Weyl/symmetric
type. It is also not necessarily of the Wick/normal type
\cite{berezin74,berezin75,berezin76,kar96,bordewald97,dollyasha01,neumaier03,
donin03}.
\item
We let the two geometric input data \mb{\mult^{ij}} and \mb{\nabla} depend on
Planck's constant \mb{\hbar}. A priori, such a generalization is a non-trivial
matter, since Planck's constant \mb{\hbar} is charged with respect to the  pertinent
resolution degree, cf.~\eq{fedosovdeg}. Moreover, for \mb{\hbar}-dependent
structures, there is no analogue of Darboux' theorem at our disposal, cf.\ discussion below \eq{ijkjacid}.
\end{enumerate}

Certain aspects of each of the three generalizations have appeared in the
literature, but never simultaneously. Our main point is that the three
generalizations taken together constitute a~natural habitat for Fedosov's
original variables \mb{x^{i}}, \mb{c^{i} \equiv dx^{i}}, \mb{y^{i}} and
\mb{\hbar}, without destroying the simplicity of his construction.

The paper starts with listing the basic setup and assumptions in
Sections~\ref{secbasicset}--\ref{secomega}. Then follows an introduction of
the relevant tools: An algebra \mb{\cA} of covariant tensors in
Sections~\ref{sectensor}--\ref{secalg}, the Fedosov resolution degree
in Section~\ref{secdeg}, the \mb{\circ} product in
Sections~\ref{seccirc}--\ref{seccirccom}, the Koszul--Tate dif\/ferential
\mb{\delta} and its cohomology in Sections~\ref{secdelta}--\ref{secdeltam1}.
Thereafter is given a discussion of Riemann curvature tensors on
supermanifolds in Sections~\ref{secriemanncurv}--\ref{secrcurv}. The f\/lat
\mb{D} connection, whose existence is guaranteed by Fedosov's 1st theorem
(Theorem~\ref{fedosov1theorem}), is discussed in
Sections~\ref{secdconn}--\ref{secflatdconn}. Fedosov's 2nd theorem
(Theorem~\ref{fedosov2theorem}), which establishes an algebra isomorphism
between symbols and horizontal zero-forms, is discussed in
Section~\ref{sechorizontalsec}. {}Finally, the \mb{*} product is constructed
in the last Section~\ref{secstarproduct}.

We shall focus on the existence of the construction and skip the important
question of uniqueness/ambiguity/equivalence/classif\/ication of \mb{*} products
for brevity. In detail, this paper is an elaboration of the material covered on
p.~138--147 in Fedosov's book~\cite{f96} subjected to the above generalizations
1--3. In particular, what we refer to as Fedosov's 1st and 2nd theorems
(Theorems~\ref{fedosov1theorem} and \ref{fedosov2theorem}) are
generalizations of Theorems~5.2.2 and 5.2.4 in \cite{f96}, respectively.
Some other references that deal with Fedosov's original construction are
\cite{pingxu98,gelretshu98,kravchenko98,farkas00,vaisman02,gadolmotos05}.
Historically, De Wilde and Lecomte were the f\/irst to prove the existence of
an associative \mb{*} product on every symplectic manifold \cite{wildecomte83}.
The same question for Poisson manifolds (which may be degenerate) was proved by
Kontsevich using ideas from string theory \cite{kontsevich9703}. Cattaneo et~al.\
gave an explicit construction in the Poisson case that merges Kontsevich's
local formula with Fedosov's f\/lat \mb{D} connection approach
\cite{catfeltom02a,catfeltom02b}.

{\em General remarks about notation.} We shall work with smooth
f\/inite-dimensional supermanifolds in terms of charts and atlases (as opposed
to, e.g.,  ringed spaces or functors of points). We shall sometimes make use of
local formulas, because these often provide the most transparent def\/initions
of sign conventions on supermanifolds. We stress that all formulas in this
paper, if not written in manifestly invariant form, hold with respect to  any coordinate
system, and they transform covariantly under general coordinate
transformations. Adjectives from supermathematics such as ``graded'',
``super'', etc., are implicitly implied. We will also follow commonly accepted
superconventions, such as, Grassmann-parity \mb{\eps} is only def\/ined modulo
\mb{2}, and ``nilpotent'' means ``nilpotent of order \mb{2}''. The sign
conventions are such that two exterior forms \mb{\xi} and \mb{\eta}, of
Grassmann-parity \mb{\eps_{\xi}}, \mb{\eps_{\eta}} and exterior form degree
\mb{p_{\xi}}, \mb{p_{\eta}}, respectively, commute in the following graded
sense
\[
\eta \wedge \xi
 = (-1)^{\eps_{\xi}\eps_{\eta}+p_{\xi}p_{\eta}}\xi\wedge\eta
\]
inside the exterior algebra. We will often not write exterior wedges
``\mb{\wedge}'' nor tensor multiplications ``\mb{\otimes}'' explicitly.
Covariant and exterior derivatives will always be from the left, while partial
derivatives can be from either left or right. We shall sometimes use round
parenthesis ``\mb{(~)}'' to indicate how far derivatives act, see e.g.,
\eqs{classicalbc0}{beq}.

\section{Basic settings and assumptions}
\label{secbasicset}

Consider a manifold \mb{M} with local coordinates \mb{x^{i}} of
Grassmann-parity \mb{\eps_{i}\!\equiv\eps(x^{i})}. The classical limit
\begin{gather}
f*g = f g +\frac{i\hbar}{2} \bigg(f\papar{x^{j}}\bigg) \mult^{jk} \bigg(\papal{x^{k}}g\bigg)
+ {\cal O}(\hbar^{2})
 ,\qquad f,g \in C^{\infty}(M)[[\hbar]] \equiv \cA^{}_{00}  ,
\label{classicalbc0}
\end{gather}
of the sought-for associative \mb{*} multiplication is prescribed by a
Grassmann-even contravariant \mb{(2,0)} tensor
\begin{gather*}
\mult = \mult^{ij}\partial^{\ell}_{j}\otimes\partial^{\ell}_{i}
 ,\qquad \eps(\mult) = 0 ,\qquad \eps(\partial^{\ell}_{i}) = \eps_{i} .
\end{gather*}
(The letter ``\mb{\mult}'' is a mnemonic for the word ``multiplication''.)
The tensor
\begin{gather*}
\mult^{ij} = \mult^{ij}(x;\hbar)
\end{gather*}
can in principle be a formal power series in Planck's constant \mb{\hbar}.
(The quantum corrections to \mb{\mult^{ij}} do not enter actively into the
classical boundary condition \e{classicalbc0}, but they will nevertheless
af\/fect the Fedosov implementation of the \mb{*} multiplication at higher orders
in \mb{\hbar}, as we shall see in \eq{circproduct} below.) Usually one demands
\cite{bffls78} that the classical unit function \mb{{\bf 1}\in C^{\infty}(M)}
serves as a unit for the full quantum algebra \mb{(\cA^{}_{00},+,*)}:
\begin{gather}
f*{\bf 1} = {\bf 1}*f = f ,\qquad
f \in C^{\infty}(M)[[\hbar]] \equiv \cA^{}_{00} .\label{starunit}
\end{gather}
(The notation \mb{\cA^{}_{00}} will be explained in Section~\ref{sectensor}.)
Let \mb{\mult^{T}} denote the transposed tensor,
\begin{gather*}
(\mult^{T})^{ij} := (-1)^{\eps_{i}\eps_{j}}\mult^{ji} .
\end{gather*}
It will be necessary to assume that the antisymmetric part
\begin{gather}
\omega^{ij} := \Hf \big(\mult^{ij}-(\mult^{T})^{ij}\big)
 = -(-1)^{\eps_{i}\eps_{j}}\omega^{ji}
\label{omegaupdef}
\end{gather}
of the tensor \mb{\mult^{ij}} is non-degenerate, i.e.,  that there exists an
inverse matrix \mb{\omega_{ij}} such that
\begin{gather*}
\omega_{ij}\omega^{jk} = \delta_{i}^{k} .
\end{gather*}
(Note that this does not necessarily imply that \mb{\mult^{ij}} itself has to
be non-degenerate.) Next, let there be given a torsion-free connection
\mb{\nabla:\Gamma(TM)\times\Gamma(TM)[[\hbar]]\to\Gamma(TM)[[\hbar]]} that
preserves the \mb{\mult}-tensor
\begin{gather}
(\nabla\mult) = 0 .\label{nablapreservm}
\end{gather}
In local coordinates, the covariant derivative \mb{\nabla_{X}} along a vector
f\/ield \mb{X=X^{i}\partial^{\ell}_{i}} is given as \cite{b97}
\begin{gather}
\nabla_{X} = X^{i} \nabla_{i} ,\qquad
\nabla_{i} = \papal{x^{i}}
+ \partial^{r}_{k} \; \Gamma^{k}{}_{ij} \larrow{dx^{j}} .\label{nabladef1}
\end{gather}
Equivalently, \mb{\nabla=d+\Gamma=dx^{i}\otimes\nabla_{i}
: \Gamma(TM) \to \Gamma(T^{*}M \otimes TM)}, where
\begin{gather}
\Gamma = dx^{i} \otimes \partial^{r}_{k}\; \Gamma^{k}{}_{ij} \larrow{dx^{j}} .
\label{nabladef2}
\end{gather}
Here \mb{\partial^{r}_{i}\!\equiv\!(-1)^{\eps_{i}}\partial^{\ell}_{i}} are not
usual partial derivatives. In particular, they do not act on the Christof\/fel
symbols \mb{\Gamma^{k}{}_{ij}} in \eqs{nabladef1}{nabladef2}. Rather they are a
dual basis to the one-forms \mb{\larrow{dx^{i}}}:
\begin{gather*}
\larrow{dx^{i}}(\partial^{r}_{j}) = \delta^{i}_{j} .
\end{gather*}
Phrased dif\/ferently, the \mb{\partial^{r}_{i}} are merely bookkeeping devices,
that transform as right partial derivatives under general coordinate
transformations. (To be able to distinguish them from true partial derivatives,
the dif\/ferentiation variable \mb{x^{i}} on a true partial derivative
\mb{\partial/\partial x^{i}} is written explicitly.)
The assumption \e{nablapreservm} reads in local coordinates
\begin{gather*}
0 = (\nabla_{i}\mult^{jk})
 = \bigg(\papal{x^{i}}\mult^{jk}\bigg)+\Gamma_{i}{}^{j}{}_{n} \mult^{nk}
+(-1)^{\eps_{j}(\eps_{k}+\eps_{n})}\Gamma_{i}{}^{k}{}_{n} \mult^{jn} ,
\end{gather*}
where we have introduced a reordered Christof\/fel symbol
\begin{gather*}
\Gamma_{i}{}^{k}{}_{j} := (-1)^{\eps_{i}\eps_{k}}\Gamma^{k}{}_{ij} .
\end{gather*}
Note that the connection \mb{\nabla} will also preserve the transposed tensor
\mb{\mult^{T}}, and therefore, by linearity, the antisymmetric part
\mb{\omega^{ij}}. We shall later explain why it is crucial that the connection~\mb{\nabla} is torsion-free, see comment after \eq{rzerocurvature}.
The Christof\/fel symbols \mb{\Gamma^{k}{}_{ij} = \Gamma^{k}{}_{ij}(x;\hbar)}
for the connection \mb{\nabla} is allowed to be a formal power series in
Planck's constant \mb{\hbar}. {}Finally, we mention that one traditionally
imposes a reality/hermiticity condition on the connection \mb{\nabla}, the
multiplicative structure \mb{\mult^{ij}} and the \mb{*} product.
However, we shall for simplicity skip the reality/hermiticity condition in
this paper.

\section[The two-form $\omega$ is symplectic]{The two-form $\boldsymbol{\omega}$ is symplectic}
\label{secomega}

The inverse matrix \mb{\omega_{ij}} with lower indices has the following
graded skewsymmetry
\begin{gather}
\omega_{ij}  = (-1)^{(\eps_{i}+1)(\eps_{j}+1)}\omega_{ji} ,
\label{omegaskewsym}
\end{gather}
cf.\  \eq{omegaupdef}. That \eq{omegaskewsym} should be counted as an
skewsymmetry (as opposed to a symmetry) is perhaps easiest to see by def\/ining
the slightly modif\/ied matrix
\begin{gather*}
\tilde{\omega}_{ij} := \omega_{ij}(-1)^{\eps_{j}} .
\end{gather*}
Note that the two matrices \mb{\omega_{ij}} and \mb{\tilde{\omega}_{ij}} are
identical for bosonic manifolds. Then the \eq{omegaskewsym} translates into
the more familiar type of graded antisymmetry,
\begin{gather}
\tilde{\omega}_{ij}  = -(-1)^{\eps_{i}\eps_{j}}\tilde{\omega}_{ji} .
\label{tildeomegaskewsym}
\end{gather}
The skewsymmetry means that the inverse matrix can be viewed as a two-form
\begin{gather}
\omega := \Hf c^{i}\omega_{ij}c^{j} = -\Hf c^{j}c^{i}\tilde{\omega}_{ij}
 \in \Gamma\left(\bigwedge{}^{2}(T^{*}M)\right)[[\hbar]] \equiv \cA_{20} .
\label{omegatwoform}
\end{gather}
Here
\begin{gather*}
c^{i} \equiv dx^{i}
\end{gather*}
is the usual basis of one-forms, and
\begin{gather*}
d := c^{i}\papal{x^{i}}
\end{gather*}
denotes the de Rham exterior derivative on \mb{M}. It follows from assumption
\e{nablapreservm} that the connection \mb{\nabla} preserves the
two-form
\begin{gather}
0 = (\nabla_{i}\tilde{\omega}_{jk})
 = \bigg(\papal{x^{i}}\tilde{\omega}_{jk}\bigg)
-\left((-1)^{\eps_{i}\eps_{j}}\Gamma_{j,ik}
-(-1)^{\eps_{j}\eps_{k}}(j\leftrightarrow k)\right) ,
\label{nablapreservomdown}
\end{gather}
where the lowered Christof\/fel symbol \mb{\Gamma_{k,ij}} is def\/ined as
\begin{gather*}
\Gamma_{k,ij} := \omega_{kn}\Gamma^{n}{}_{ij}(-1)^{\eps_{j}} .
\end{gather*}
The two-form \mb{\omega} is closed
\begin{gather}
(d\omega) = 0 ,
\label{omegaclosed}
\end{gather}
or equivalently,
\begin{gather*}
\sum_{{\rm cycl.}~i,j,k}(-1)^{\eps_{i}\eps_{k}}
\bigg(\papal{x^{i}}\tilde{\omega}_{jk}\bigg) = 0 .
\end{gather*}
The closeness relation \e{omegaclosed} is not an extra assumption. It follows
from \eq{nablapreservomdown}, because the connection \mb{\nabla} is
torsion-free,
\begin{gather}
T^{k}{}_{ij}
 := \Gamma^{k}{}_{ij}+(-1)^{(\eps_{i}+1)(\eps_{j}+1)}(i\leftrightarrow j)
 = 0 ,\label{torsionfreeup}
\end{gather}
or equivalently, in terms of the lowered Christof\/fel symbol,
\begin{gather*}
\Gamma_{k,ij} = (-1)^{\eps_{i}\eps_{j}}(i\leftrightarrow j) .
\end{gather*}
It is practical to call a non-degenerate closed two-form \mb{\omega_{ij}} a
{\em symplectic} structure, even if it depends on Planck's constant \mb{\hbar}.
The inverse structure, i.e.,  the corresponding Poisson structure~\mb{\omega^{ij}}
satisf\/ies the Jacobi identity
\begin{gather}
\sum_{{\rm cycl.}~i,j,k}(-1)^{\eps_{i}\eps_{k}}
\omega^{in} \bigg(\papal{x^{n}}\omega^{jk}\bigg)  = 0 .\label{ijkjacid}
\end{gather}
Note that we cannot rely on Darboux' theorem, i.e.,  we will not be guaranteed a
cover of Darboux coordinate patches in which the \mb{\omega^{ij}} is constant.
The issue is that, on one hand, the symplectic structure is allowed to depend
on Planck's constant~\mb{\hbar}, but, on the other hand, we shall not allow
coordinate transformations \mb{x^{i}\to x^{\prime j}=x^{\prime j}(x)} that
depend on~\mb{\hbar}. Luckily, as we shall see, Darboux patches play no r\^ole
in the Fedosov construction. In fact, as we have mentioned before, all
formulas in this paper, if not written in manifestly invariant form, hold with respect to
any coordinate system, and they transform covariantly under general coordinate
transformations.

The classical Poisson bracket is given by the famous {\em quantum
correspondence principle} \cite{dirac25}
\begin{gather}
\{f,g\}^{}_{\cl} := \bigg(f\papar{x^{j}}\bigg)\omega_{(0)}^{jk}\bigg(\papal{x^{k}}g\bigg)
 = \lim_{\hbar\to 0}\frac{1}{i\hbar}[f\starcomma g]
 ,\qquad   f,g \in C^{\infty}(M) .
\label{classicalpb}
\end{gather}
Here
\begin{gather*}
[f\starcomma g] := f*g -(-1)^{\eps_{f}\eps_{g}} g*f
 = i\hbar\{f,g\}^{}_{\cl}+ {\cal O}(\hbar^{2})
 ,\qquad f,g\in C^{\infty}(M)[[\hbar]] \equiv \cA^{}_{00} ,
\end{gather*}
is the \mb{*} commutator, and
\begin{gather*}
\omega_{(0)}^{ij} := \lim_{\hbar\to 0}\omega^{ij} .
\end{gather*}
It is easy to show that every symplectic manifold \mb{(M;\omega)} has a
torsion-free \mb{\omega}-preserving connection \mb{\nabla}, see Section~2.5 in~\cite{f96} for the bosonic case. However, it is not true that every manifold
\mb{(M;\mult)} with a multiplicative structure \mb{\mult^{ij}} supports a
torsion-free \mb{\mult}-preserving connection \mb{\nabla}, cf.\  assumption
\eq{nablapreservm}. The symmetric part
\begin{gather*}
g^{ij} := \Hf (\mult^{ij}+(\mult^{T})^{ij}) = (-1)^{\eps_{i}\eps_{j}}g^{ji}
\end{gather*}
of the tensor \mb{\mult^{ij}} needs to be compatible with the symplectic
structure \mb{\omega} in a certain sense. In the special case where
\mb{g^{ij}=0}, we return to the usual Fedosov quantization
\mb{\mult^{ij}=\omega^{ij}}, which corresponds to a Weyl/symmetric type \mb{*}
product. In the generic case where \mb{g^{ij}} has full rank, there will
exist an inverse matrix \mb{g_{ij}}, which constitute a (pseudo) Riemannian
metric, and there will hence exist a corresponding unique Levi-Civita
connection \mb{\nabla^{LC}}. In this non-degenerate case, the necessary and
suf\/f\/icient conditions are \mb{\nabla=\nabla^{LC}} and
\mb{(\nabla^{LC}_{i}\omega^{jk})=0}. This is for instance satisf\/ied for
(pseudo) K\"ahler manifolds \mb{(M;\omega;g)}, cf.\
\cite{berezin74,berezin75,berezin76,kar96,bordewald97,dollyasha01,neumaier03,
donin03}.

\section{Covariant tensors}
\label{sectensor}

Let
\begin{gather}
\cA_{\mm\nn} := \Gamma\left(\bigwedge{}^{\mm}(T^{*}M)
\otimes\bigveee^{\nn}(T^{*}M)\right)[[\hbar]]
\label{anmalgebra}
\end{gather}
be the vector space of \mb{(0,\mm + \nn)}-tensors
\mb{a_{i_{1}\dots i_{\mm}j_{1}\dots j_{\nn}}(x;\hbar)}
that are antisymmetric with respect to the  f\/irst \mb{\mm} indices \mb{i_{1}\ldots i_{\mm}}
and symmetric with respect to the  last \mb{\nn} indices \mb{j_{1}\ldots j_{\nn}}. Phrased
dif\/ferently, they are \mb{\mm}-form valued symmetric \mb{(0,\nn)}-tensors.
As usual, it is practical to introduce a coordinate-free notation
\begin{gather}
a_{\mm\nn}(x,c;y;\hbar)
 = \frac{1}{\mm!\nn!}c^{i_{\mm}}\wedge\cdots\wedge c^{i_{1}}
a_{i_{1}\dots i_{\mm}j_{1}\dots j_{\nn}}(x;\hbar)\otimes y^{j_{\nn}}
\veee\cdots\veee y^{j_{1}} .
\label{coordinatefree}
\end{gather}
Here the Fedosov variables \mb{y^{i}} and the ``\mb{\veee}'' symbol are
symmetric counterparts to the one-form basis \mb{c^{i}\equiv dx^{i}} and the
``\mb{\wedge}'' symbol, i.e.,
\begin{gather*}
y^{j}\veee y^{i} = (-1)^{\eps_{i}\eps_{j}} y^{i}\veee y^{j} .
\end{gather*}
We will be interested in covariant derivatives
\mb{\nabla_{i}a_{i_{1}\dots i_{\mm}j_{1}\dots j_{\nn}}} of the above tensors.
The covariant derivative \mb{\nabla_{i}} from \eq{nabladef1} can be implemented
on coordinate-free objects \e{coordinatefree} by the following linear
dif\/ferential operator
\begin{gather}
\nabla_{i} = \papal{x^{i}}-\Gamma_{i}{}^{k}{}_{j}c^{j}\papal{c^{k}}
-\Gamma_{i}{}^{k}{}_{j}y^{j}\papal{y^{k}} .
\label{nablaidef}
\end{gather}
If both the numbers of antisymmetric and symmetric indices are non-zero
\mb{\mm\neq 0 \wedge \nn\neq 0}, i.e.,  if the tensor
\mb{a_{i_{1}\dots i_{\mm}j_{1}\dots j_{\nn}}} has mixed symmetry, the
covariant derivative
\mb{\nabla_{i}a_{i_{1}\dots i_{\mm}j_{1}\dots j_{\nn}}} will not belong to
any of the \mb{\cA_{\bullet\bullet}} spaces \e{anmalgebra}. We repair this by
antisymmetrizing with respect to the  indices \mb{i,i_{1},\ldots,i_{\mm}}. Such
antisymmetrization can be implemented on coordinate-free objects~\e{coordinatefree} with the help of the  following one-form valued Grassmann-even
dif\/ferential operator
\begin{gather}
\nabla := c^{i}\nabla_{i} = d-c^{i}\Gamma_{i}{}^{k}{}_{j}y^{j}\papal{y^{k}} ,
\label{nabladef}
\end{gather}
where we have followed common practice, and given the dif\/ferential operator
\e{nabladef} the same name as the connection. In the second equality of
\eq{nabladef} we have used that the connection is torsion-free. (References~\cite{karsch01,kar03} consider a hybrid model where torsion is allowed
in the \mb{y}-sector but not in the \mb{c}-sector in such a way that
\eq{nabladef} remains valid.) Since the \mb{\nabla} operator is a~f\/irst-order
dif\/ferential operator, i.e.,
\begin{gather}
\nabla(ab) = (\nabla a)b+(-1)^{p_{a}}a(\nabla b) , \label{nablalinear}
\end{gather}
where \mb{a} and \mb{b} are two coordinate-free objects \e{coordinatefree}, it
is customary to refer to \mb{\nabla} as a {\em linear} connection. (The order
of the exterior factor \mb{\bigwedge{}^{\mm}(T^{*}M)} and the symmetric factor
\mb{\bigveee^{\nn}(T^{*}M)} in expression \e{anmalgebra} is opposite the
standard convention to ease the use of covariant derivatives~\mb{\nabla} that
acts from the {\em left}.)

\begin{table}[t]
\caption{Parities and gradings.}\label{gradtable}
\begin{center}
\begin{tabular}{|l|c||c|c|c|}  \hline
&Grading $\rightarrow$&Grassmann&Form&Fedosov \\
&&parity&degree&degree\\ \hline
$\downarrow$ Variable& $\downarrow$ Symbol $\rightarrow$&$\eps$&$p$&$\deg$\\
\hline\hline
Coordinates&$x^{i}$&$\eps_{i}$&$0$&$0$ \\ \hline
$1$-form basis&$c^{i}\!\equiv\!dx^{i}$&$\eps_{i}$&$1$&$0$ \\ \hline
{}Fedosov coordinates&$y^{i}$&$\eps_{i}$&$0$&$1$ \\ \hline
Planck's constant&$\hbar$&$0$&$0$&$2$  \\ \hline\hline
Multiplicative structure&$\mult^{ij}$&$\eps_{i}\!+\!\eps_{j}$&$0$
&$\even,\geq 0$ \\ \hline
Christof\/fel symbol&$\Gamma^{k}{}_{ij}$&$\eps_{i}\!+\!\eps_{j}\!+\!\eps_{k}$&$0$
&$\even,\geq 0$ \\ \hline
Covariant derivative&$\nabla_{i}$&$\eps_{i}$&$0$&$\even,\geq 0$ \\ \hline
$1$-form valued connection&$\nabla\!\equiv\!c^{i}\nabla_{i}$&$0$&$1$
&$\even,\geq 0$ \\ \hline
de Rham exterior derivative &$d$&$0$&$1$&$0$ \\ \hline
Koszul--Tate dif\/ferential&$\delta\!=\!\{\omegaone,\cdot\}$&$0$&$1$&$-1$\\ \hline
Contraction&$\delta^{*}\!=\!y^{j} i(\partial^{\ell}_{j})$&$0$&$-1$&$1$ \\
\hline
Koszul--Tate Hamiltonian&$\omegaone\!\equiv\!c^{i}\omega_{ij}y^{j}$&$0$&$1$
&$\odd,\geq 1$ \\ \hline
Deformation $1$-form&$r$&$0$&$1$&$\geq 0$ \\ \hline
Hamiltonian curvature $2$-form&$\cR$&$0$&$2$&$\even,\geq 2$ \\ \hline
\end{tabular}
\end{center}
\end{table}

\section[The $\cA$ algebra]{The $\boldsymbol{\cA}$ algebra}
\label{secalg}

The direct sum
\begin{gather}
\cA := \bigoplus_{\mm,\nn\geq 0}\cA_{\mm\nn}
\cong C^{\infty}(E)[[\hbar]]
\label{aalgebra}
\end{gather}
of the \mb{\cA_{\mm\nn}} spaces \e{anmalgebra} is an algebra with
multiplication given by the tensor multiplication. It is both associative and
commutative. As indicated in \eq{aalgebra}, the elements
\begin{gather*}
a = \oplus_{\mm,\nn\geq 0}a_{\mm\nn} \in \cA
 ,\qquad   a_{\mm\nn} \in \cA_{\mm\nn} ,
\end{gather*}
can be viewed as quantum functions \mb{a=a(x,c;y;\hbar)} on the Whitney sum
\begin{gather*}
E := \Pi TM\oplus TM ,
\end{gather*}
where \mb{c^{i}} are identif\/ied with the parity-inverted f\/iber coordinates for
the f\/iber-wise parity-inverted tangent bundle~\mb{\Pi TM}, and~\mb{y^{j}} are
identif\/ied with the f\/iber coordinates for the tangent bundle~\mb{TM}. The
word ``quantum function'' just means that it is a formal power series in
Planck's constant~\mb{\hbar}.

\section{The Fedosov resolution degree}
\label{secdeg}

The {\em Fedosov degree} ``\mb{\deg}'' is a (non-negative) integer grading of
the \mb{\cA} algebra def\/ined as
\begin{gather}
\deg(y^{i}) = 1 ,\qquad   \deg(\hbar) = 2 ,
\label{fedosovdeg}
\end{gather}
and zero for the two other remaining variables \mb{x^{j}} and \mb{c^{k}}, cf.\
Table~\ref{gradtable}. The Fedosov degree will play the r\^ole of resolution
degree in the sense of homological perturbation theory
\cite{fravil75,batvil77,frafra78,batfra83a,batfra83b,fhst89,hentei92}. We
shall therefore often organize the algebra according to this grading:
\begin{gather*}
\cA = \bigoplus_{\nn\geq 0}\cA_{(\nn)} ,\qquad
\cA_{(\nn)} := \{a\in\cA \mid \deg(a)=\nn\} .
\end{gather*}
Similarly, one may write the algebra element
\begin{gather*}
a = \oplus_{\nn\geq 0}a_{(\nn)} \in \cA
 ,\qquad   a_{(\nn)} \in \cA_{(\nn)} ,
\end{gather*}
as a direct sum of elements \mb{a_{(\nn)}} of def\/inite Fedosov degree
\mb{\deg(a_{(\nn)})=\nn}.

\section[The $\circ$ product]{The $\boldsymbol{\circ}$ product}
\label{seccirc}

One now builds an associative \mb{\circ} product \mb{\cA\times\cA\to\cA} on
the \mb{\cA} algebra as a Moyal product \cite{groen46,moyal49} in the
\mb{y}-variables,
\begin{gather}
a\circ b :=
\left(a\exp\left[\papar{y^{j}}\frac{i\hbar}{2}\mult^{jk}\papal{y^{k}}\right]b\right)
 ,\qquad   a,b\in \cA .
\label{circproduct}
\end{gather}
The \mb{\circ} product is associative, because the \mb{\mult^{jk}}-tensor is
independent of \mb{y}-variables. (The \mb{y}-variables have been interpreted
by Grigoriev and Lyakhovich \cite{grilya01,batgrilya01} as conversion
variables for the conversion of second-class constraints into f\/irst-class
\cite{bf87,bff89,bt91,fl94}.) The \mb{\circ} product ``preserves'' the
following gradings
\begin{gather*}
\deg(a\circ b) \geq \deg(a)+\deg(b) ,\\ 
p(a\circ b) = p(a)+p(b) , \\ 
\eps(a\circ b) = \eps(a)+\eps(b) .
\end{gather*}
The connection \mb{\nabla} respects the \mb{\circ} product,
\begin{gather}
\nabla(a\circ b) = (\nabla a)\circ b+(-1)^{p_{a}}a\circ(\nabla b)
 ,\qquad   a,b\in \cA ,
\label{nablarespectcirc}
\end{gather}
as a consequence of the assumption \e{nablapreservm}.

\section{The Poisson bracket}
\label{secpb}

It is useful to def\/ine a Poisson bracket as
\begin{gather}
\{a,b\} := \bigg(a\papar{y^{i}}\bigg)\omega^{ij}\bigg(\papal{y^{j}}b\bigg)
 ,\qquad   a,b\in \cA .\label{poissonbracket}
\end{gather}
(It should not be confused with the classical Poisson bracket \e{classicalpb}.)
The Poisson bracket \e{poissonbracket} ``preserves'' the following gradings
\begin{gather*}
0 \leq \deg(\{a, b\}) \geq \deg(a)+\deg(b)-2 , \\ 
p(\{a, b\}) = p(a)+p(b) , \\ 
\eps(\{a, b\}) = \eps(a)+\eps(b) .
\end{gather*}
The connection \mb{\nabla} respects the Poisson bracket:
\begin{gather*}
\nabla\{a, b\} = \{\nabla a, b\}+(-1)^{p_{a}}\{a,\nabla b\}
 ,\qquad   a,b\in \cA .
\end{gather*}

\section[The $\circ$ commutator]{The $\boldsymbol{\circ}$ commutator}
\label{seccirccom}

The \mb{\circ} commutator is def\/ined as
\begin{gather*}
[a\ocomma b] := a\circ b -(-1)^{\eps_{a}\eps_{b}+p_{a}p_{b}} b\circ a
 = i\hbar\{a,b\}+ {\cal O}(\hbar^{2}) ,\qquad   a,b\in \cA .
\end{gather*}
Note the following useful observations:
\begin{itemize}\itemsep=0pt
\item
Each term in the \mb{\circ} commutator \mb{[a\ocomma b]} contains at least one
power of \mb{\hbar}, so one may consider the fraction
\mb{\frac{1}{i\hbar}[a\ocomma b]} without introducing negative powers of
\mb{\hbar}.
\item
The \mb{\circ} commutator \mb{\frac{1}{i\hbar}[a\ocomma b]} and the Poisson
bracket \mb{\{a,b\}} are equal, if one of the entries \mb{a} or \mb{b} only
contains terms with less than three \mb{y}'s.
\end{itemize}

The \mb{\circ} commutator may be expanded according to the Fedosov degree:
\begin{gather*}
[a\ocomma b] = \sum_{\nn\geq 0}[a\ocomma b]_{(\nn)} ,
\end{gather*}
where
\begin{gather*}
\deg([a\ocomma b]_{(\nn)}) = \deg(a)+\deg(b)+\nn .
\end{gather*}

\section[The Koszul-Tate differential $\delta$]{The Koszul--Tate dif\/ferential $\boldsymbol{\delta}$}
\label{secdelta}

The Koszul--Tate dif\/ferential is def\/ined as
\begin{gather}
\delta := c^{i}\papal{y^{i}} = \{\omegaone,\cdot\}
 ,\qquad    \deg(\delta) = -1 .
\label{deltadef}
\end{gather}
In the second equality in \eq{deltadef} we have indicated that the Koszul--Tate
dif\/ferential is an inner derivation in the algebra
\mb{(\cA,+,\frac{1}{i\hbar}[ \cdot\ocomma\cdot ])} with generator
\begin{gather*}
\omegaone := c^{i}\omega_{ij}y^{j}
 = y^{j}c^{i}\tilde{\omega}_{ij} \in \cA_{11}
 ,\qquad    \deg(\omegaone) \geq 1 ,
\end{gather*}
which we shall refer to as the Hamiltonian for \mb{\delta}. The Koszul--Tate
dif\/ferential \mb{\delta} will serve as the leading term in a resolution
expansion of a deformed connection \mb{D}, see Section~\ref{secdconn}
\cite{fravil75,batvil77,frafra78,batfra83a,batfra83b,fhst89,hentei92}. Since
the \mb{\omega_{ij}} tensor is covariantly preserved, cf.\
\eq{nablapreservomdown}, it follows immediately that
\begin{gather*}
(\nabla\omegaone) = -y^{j}c^{i}c^{k}(\nabla_{k}\tilde{\omega}_{ij}) = 0
\end{gather*}
even without using the skewsymmetry \e{omegaskewsym} or the antisymmetry
\e{tildeomegaskewsym} (or the torsion-free condition for that matter). As a
corollary,
\begin{gather*}
[\nabla,\delta] = 0 .
\end{gather*}
The \mb{\delta}-dif\/ferential is nilpotent
\begin{gather}
\delta^{2} \equiv \Hf[\delta,\delta] = 0 ,
\label{deltanilp}
\end{gather}
and it respects the \mb{\circ} product
\begin{gather}
\delta(a\circ b) = (\delta a)\circ b+(-1)^{p_{a}}a\circ(\delta b) .
\label{deltarespectcirc}
\end{gather}

\section[The Poincar\'e lemma and the homotopy operator
$\delta^{-1}$]{The Poincar\'e lemma and the homotopy operator
$\boldsymbol{\delta^{-1}}$}
\label{secdeltam1}

There exists a version of the Poincar\'e lemma where the r\^ole of the de Rham
exterior derivative \mb{d \equiv c^{i}\partial^{\ell}/\partial x^{i}} has
been replaced by the Koszul--Tate dif\/ferential
\mb{\delta \equiv c^{i}\partial^{\ell}/\partial y^{i}}, or equivalently,
where the \mb{x}-coordinates are replaced by the \mb{y}-coordinates. As we
shall see below, the equation
\begin{gather}
(\delta b) = a ,\qquad a\in\cA, \label{beq}
\end{gather}
may be solved locally with respect to  an algebra element \mb{b\in\cA} whenever
\mb{a\in\cA} is a given \mb{\delta}-closed algebra element with no \mb{00}-part
\mb{a^{}_{00} = 0}. (The \mb{00}-part \mb{a^{}_{00}} is the part of the
algebra element~\mb{a} that is independent of the \mb{c}- and the
\mb{y}-variables, cf.\  Section~\ref{sectensor}.) In fact, a local solution
\mb{b} to \eq{beq} may be extended to a global solution, since the total space
is contractible in the \mb{y}-directions. We shall use this crucial fact to
guarantee the existence of global solutions to dif\/ferential equations, whose
dif\/ferential operator is a deformation of the Koszul--Tate dif\/ferential~\mb{\delta}, see Theorems~\ref{fedosov1theorem} and \ref{fedosov2theorem}. As
usual when dealing with the Poincar\'e Lemma, it is useful to consider the
inner contraction
\begin{gather*}
\delta^{*} := y^{j}\papal{c^{j}} = y^{j} i(\partial^{\ell}_{j}) ,
\end{gather*}
which is dual to \mb{\delta} with respect to  a \mb{y\leftrightarrow c} permutation. The
commutator
\begin{gather*}
[\delta,\delta^{*}] = c^{i}\papal{c^{i}}+y^{i}\papal{y^{i}}
\end{gather*}
is a Euler vector f\/ield that counts the number of \mb{c}'s and \mb{y}'s.
The homotopy operator \mb{\delta^{-1}} is def\/ined as
\begin{gather*}
\forall \, a\in \cA_{\mm\nn}: \delta^{-1}a := \left\{\begin{array}{lcl}
\frac{1}{\mm+\nn}(\delta^{*}a)&{\rm for}&(\mm,\nn)\neq(0,0),\cr
0&{\rm for}&(\mm,\nn)=(0,0) ,\end{array}\right.
\end{gather*}
and extended by linearity to the whole algebra \mb{\cA}. The homotopy operator
\mb{\delta^{-1}} is {\em not} a f\/irst-order dif\/ferential operator, in contrast
to \mb{\delta^{*}}. One easily obtains the following version of the Poincar\'e
lemma.

\begin{lemma}[Poincar\'e lemma]\label{poincarelemma}
There is only non-trivial \mb{\delta}-cohomology in the \mb{00}-sector with
neither \mb{c}'s nor \mb{y}'s. A more refined statement is the following:
{}For all \mb{\delta}-closed elements \mb{a \in \cA} with no \mb{00}-part,
there exists a unique \mb{\delta^{*}}-closed element \mb{b \in \cA} with no
\mb{00}-part, such that \mb{a=(\delta b)}. Or equivalently, in symbols:
\begin{gather*}
\forall \, a\in \cA:
\left\{\begin{array}{c} (\delta a) = 0 \cr
a^{}_{00} = 0 \end{array}\right\}
  \  \Rightarrow    \ \exists !\,  b\in \cA:
\left\{\begin{array}{c}a = (\delta b)\cr
(\delta^{*}b) = 0\cr b^{}_{00} = 0  \end{array}\right\} .
\end{gather*}
The unique element \mb{b} is given by the homotopy operator \mb{\delta^{-1}a}.
\end{lemma}

\section{The Riemann curvature}
\label{secriemanncurv}

See   \cite{b97} and \cite{lavrov04} for related discussions. The Riemann
curvature \mb{R} is def\/ined as (half) the commutator of the \mb{\nabla}
connection \e{nabladef2},
\begin{gather}
 R = \Hf [\nabla \stackrel{\wedge}{,}\nabla]
 = -\Hf dx^{j} \wedge dx^{i} \otimes [\nabla_{i},\nabla_{j}]
 = -\Hf dx^{j} \wedge dx^{i} \otimes \partial^{r}_{n} \; R^{n}{}_{ijk}
\larrow{dx^{k}} ,\label{riemanntensordef1}  \\
R^{n}{}_{ijk} = \larrow{dx^{n}}
\left([\nabla_{i},\nabla_{j}]\partial^{r}_{k}\right)  = (-1)^{\eps_{n}\eps_{i}}\bigg(\papal{x^{i}}\Gamma^{n}{}_{jk}\bigg)
+\Gamma^{n}{}_{im}\Gamma^{m}{}_{jk}
-(-1)^{\eps_{i}\eps_{j}}(i\leftrightarrow j) ,\label{riemanntensor1}
\end{gather}
where it is implicitly understood in the second equality of
\eq{riemanntensordef1} that \mb{\nabla_{i}} does not act on \mb{dx^{j}}.
If \mb{\nabla_{i}} is supposed to act on \mb{dx^{j}}, one should include a
torsion-term
\begin{gather}
dx^{i}\wedge [\nabla_{i},dx^{k}] \nabla_{k}
 = -dx^{i}\wedge\Gamma_{i}{}^{k}{}_{j} dx^{j} \nabla_{k}
 = \frac{(-1)^{\eps_{i}}}{2}T^{k}{}_{ij} dx^{j}\wedge dx^{i} \nabla_{k} ,
\label{additionaltorsionterm}
\end{gather}
cf.\  \eq{curvasham}. Note that the order of indices in the Riemann curvature
tensor \mb{R^{n}{}_{ijk}} is non-standard. This is to minimize appearances of
Grassmann sign factors. Alternatively, the Riemann curvature tensor may be
def\/ined as
\begin{gather*}
R(X,Y)Z = \left([\nabla_{X},\nabla_{Y}]-\nabla_{[X,Y]}\right)Z
 = Y^{j}X^{i}R_{ij}{}^{n}{}_{k}Z^{k} \partial^{\ell}_{n} ,
\end{gather*}
where \mb{X=X^{i}\partial^{\ell}_{i}}, \mb{Y=Y^{j}\partial^{\ell}_{j}} and
\mb{Z=Z^{k}\partial^{\ell}_{k}} are vector f\/ields of even Grassmann-parity.
The Riemann curvature tensor \mb{R_{ij}{}^{n}{}_{k}} reads in local coordinates
\begin{gather*}
R_{ij}{}^{n}{}_{k} = (-1)^{\eps_{n}(\eps_{i}+\eps_{j})}R^{n}{}_{ijk}
 = \bigg(\papal{x^{i}}\Gamma_{j}{}^{n}{}_{k}\bigg)
+(-1)^{\eps_{j}\eps_{n}}\Gamma_{i}{}^{n}{}_{m}\Gamma^{m}{}_{jk}
-(-1)^{\eps_{i}\eps_{j}}(i\leftrightarrow j) .
\end{gather*}
It is sometimes useful to reorder the indices in the Riemann curvature tensors
as
\begin{gather}
R_{ijk}{}^{n} = ([\nabla_{i},\nabla_{j}]\partial^{\ell}_{k})^{n}
 = (-1)^{\eps_{k}(\eps_{n}+1)}R_{ij}{}^{n}{}_{k} .
\label{riemanntensor4}
\end{gather}
For a symplectic connection \mb{\nabla}, we prefer to work with a \mb{(0,4)}
Riemann tensor (as opposed to a \mb{(1,3)} tensor) by lowering the upper index
with the symplectic metric \e{omegatwoform}. In terms of Christof\/fel symbols
it is easiest to work with expression \e{riemanntensor1}:
\begin{gather}
R_{n,ijk} := \omega_{nm}R^{m}{}_{ijk}(-1)^{\eps_{k}}\nonumber\\
 \phantom{R_{n,ijk}  \:}{} = (-1)^{\eps_{i}\eps_{n}}\bigg(\papal{x^{i}}\Gamma_{n,jk}
+(-1)^{\eps_{m}(\eps_{i}+\eps_{n}+1)+\eps_{k}}
\Gamma_{m,in}\Gamma^{m}{}_{jk}\bigg)   -(-1)^{\eps_{i}\eps_{j}}(i\leftrightarrow j) .\label{riemanntensor1low}
\end{gather}
In the second equality of \eq{riemanntensor1low} we have use the symplectic
condition \e{nablapreservomdown}. If we use the symplectic condition
\e{nablapreservomdown} one more time on the f\/irst term in
\eq{riemanntensor1low}, we derive the following symmetry
\begin{gather}
R_{n,ijk}
 = (-1)^{(\eps_{i}+\eps_{j})(\eps_{k}+\eps_{n})+\eps_{k}\eps_{n}}
(k\leftrightarrow n) .
\label{riemanntensor1lowsym}
\end{gather}
This symmetry becomes clearer if we instead start from expression
\e{riemanntensor4} and def\/ine
\begin{gather*}
R_{ij,kn} := R_{ijk}{}^{m}\tilde{\omega}_{mn}
 = -(-1)^{\eps_{n}(\eps_{i}+\eps_{j}+\eps_{k})}R_{n,ijk} .
\end{gather*}
Then the symmetry \e{riemanntensor1lowsym} simply translates into a symmetry
between the third and fourth index:
\begin{gather*}
R_{ij,kn} = (-1)^{\eps_{k}\eps_{n}}(k\leftrightarrow n) .
\end{gather*}
We note that the torsion-free condition has not been used at all in this
Section~\ref{secriemanncurv}.

\section[The curvature two-form $R$]{The curvature two-form $\boldsymbol{R}$}
\label{secrcurv}

Let us now calculate the commutator of two \mb{\nabla_{i}} operators using the
realization \e{nablaidef} of the covariant derivative:
\begin{gather*}
-[\nabla_{i},\nabla_{j}]
 = R_{ij}{}^{n}{}_{k}c^{k}\papal{c^{n}}
+R_{ij}{}^{n}{}_{k}y^{k}\papal{y^{n}}
 = R_{ij}{}^{n}{}_{k}c^{k}\papal{c^{n}}
+\Hf (-1)^{\eps_{k}+\eps_{n}} \{R_{ij,kn}y^{n}y^{k}, \cdot\} .\!\!\!\!
\end{gather*}
Using realization \e{nablaidef} the curvature two-form reads
\begin{gather}
R := \nabla^{2} = c^{i}[\nabla_{i},c^{k}]\nabla_{k}
-\Hf c^{j}c^{i}[\nabla_{i},\nabla_{j}] \nonumber\\
\phantom{R \:}{}  = \frac{(-1)^{\eps_{i}}}{2}T^{k}{}_{ij}c^{j}c^{i}\nabla_{k}
+\Hf c^{k}c^{j}c^{i}R_{ijk}{}^{n}\papal{c^{n}}
+\Hf y^{k}c^{j}c^{i}R_{ijk}{}^{n}\papal{y^{n}} \nonumber\\
\phantom{R \: }{} = \Hf y^{k}c^{j}c^{i}R_{ij,kn}(-1)^{\eps_{n}}\{y^{n},\cdot\}
 = \{\cR,\cdot\} ,\label{curvasham}
\end{gather}
where
\begin{gather*}
\cR := \frac{1}{4}y^{n}y^{k}c^{j}c^{i}R_{ij,kn} \in \cA_{22}
\end{gather*}
is a Hamiltonian for the curvature two-form~\mb{R}. In this formalism the
connection \mb{\nabla_{i}} is supposed to act on \mb{c^{k}\!\equiv\!dx^{k}}, so
the torsion-term \e{additionaltorsionterm} is included in the curvature
two-form \e{curvasham}. In the fourth equality of~\eq{curvasham} we use the
torsion-free condition \e{torsionfreeup} and the f\/irst Bianchi identity for a
torsion-free connection
\begin{gather}
0 = \sum_{i,j,k~{\rm cycl.}}(-1)^{\eps_{i}\eps_{k}}R_{ijk}{}^{n}
\label{1bianchiid}
\end{gather}
to ensure that the~\mb{x}- and \mb{c}-derivative term in~\eq{curvasham}
vanishes. As we shall see below, it is vital that there is no \mb{x}- and
\mb{c}-derivative term in the \mb{R} curvature~\e{curvasham}. This is the main
reason why the connection \mb{\nabla} is assumed to be torsion-free. (See also
the comment after \eq{rzerocurvature}.) The f\/irst Bianchi identity
\e{1bianchiid} also implies that the Hamiltonian curvature two-form \mb{\cR} is
\mb{\delta}-closed:
\begin{gather*}
(\delta\cR) = \Hf c^{k}c^{j}c^{i}R_{ij,kn}y^{n}(-1)^{\eps_{n}} = 0 .
\end{gather*}
Similarly, the second Bianchi identity for a torsion-free connection
\begin{gather*}
0 = \sum_{i,j,k~{\rm cycl.}}(-1)^{\eps_{j}\eps_{k}}
\nabla_{k}R_{ij}{}^{n}{}_{m}
\end{gather*}
implies that the \mb{\nabla} operator preserves the Hamiltonian curvature
two-form~\mb{\cR}:
\begin{gather*}
(\nabla\cR) = \frac{1}{4} c^{j}c^{i}c^{k}(\nabla_{k}R_{ij,nm})
y^{m}y^{n}(-1)^{\eps_{n}+\eps_{m}} = 0 .
\end{gather*}

\section[Higher-order $D$ connection]{Higher-order $\boldsymbol{D}$ connection}
\label{secdconn}

We next deform the linear \mb{\nabla} connection from equation \e{nabladef}
into a higher-order connection \mb{D:\Gamma(TM)[[\hbar]]\times\cA\to\cA},
\begin{gather*}
D := \nabla-\delta+\frac{1}{i\hbar}[r\ocomma\cdot \,]
 = \nabla+\frac{1}{i\hbar}[r\!-\!\omegaone\ocomma\cdot \, ]
 = \sum_{\nn\geq -1}\stackrel{(\nn)}{D}
 ,\qquad \deg(\stackrel{(\nn)}{D}) = \nn ,
\end{gather*}
with the help of  a so-called {\em deformation one-form}
\begin{gather*}
r = \oplus_{\nn\geq 0}r_{(\nn)} \in \cA_{1\bullet} ,
\qquad    \deg(r_{(\nn)}) = \nn .
\end{gather*}
As we soon shall see, it is better to think of \mb{D} as a deformation of
(minus) the Koszul--Tate dif\/ferential \mb{\delta} rather than the connection
\mb{\nabla}. The word ``higher-order'' refers to that \mb{D} is {\em not}
necessarily a linear derivation of the tensor algebra \mb{(\cA,+,\otimes)},
cf.\  \eq{nablalinear}. However, it is a linear derivation of the \mb{\circ}
algebra \mb{(\cA,+,\circ)}, similar to
\eqs{nablarespectcirc}{deltarespectcirc}. The \mb{\nabla} and~\mb{D}
connections may be expanded in Fedosov degree:
\begin{gather*}
\nabla = \sum_{\nn\geq 0}\stackrel{(\nn)}{\nabla}
 ,\qquad    \deg(\stackrel{(\nn)}{\nabla}) = \nn
 ,\qquad    \stackrel{(\odd)}{\nabla} = 0 ,
\end{gather*}
and
\begin{gather}
\stackrel{(\nn)}{D} := \left\{\begin{array}{lcl}
\displaystyle \frac{1}{i\hbar}[r_{(1)}\ocomma\cdot\, ]_{(0)}-\delta&{\rm for}&\nn=-1 ,\vspace{2mm}\\
\displaystyle \stackrel{(\nn)}{\nabla}+\frac{1}{i\hbar}\sum_{k=0}^{\nn+2}[r_{(k)}
\ocomma\cdot \, ]_{(\nn+2-k)}&{\rm for}&\nn\geq 0 .\end{array}\right.
\label{dmconn}
\end{gather}
Note that the connection \mb{D} does not depend on \mb{r_{(0)}}. Also note
that it will be necessary to assume that the \mb{(1)}-sector vanishes
\begin{gather*}
r_{(1)} = 0 
\end{gather*}
to ensure that (minus) the Koszul--Tate dif\/ferential \mb{\delta} is the
sole leading term in the~\mb{D} expansion~\e{dmconn}.

\section[The $R_{D}$ curvature]{The $\boldsymbol{R_{D}}$ curvature}
\label{secrdcurv}

The curvature two-form \mb{R_{D}} for the \mb{D} connection is
\begin{gather}
R_{D} := D^{2}
 = \left(\nabla+\frac{1}{i\hbar}[r\!-\!\omegaone\ocomma\cdot\, ]\right)^{2} \nonumber\\
 \phantom{R_{D} \: }{} = \nabla^{2}+\frac{1}{i\hbar}[\nabla(r\!-\!\omegaone)\ocomma\cdot \, ]
+\frac{1}{(i\hbar)^2}[r\!-\!\omegaone\ocomma [r\!-\!\omegaone\ocomma\cdot \, ]]
 = \frac{1}{i\hbar}[\cR_{D}\ocomma\cdot\, ] ,\label{rdcurv}
\end{gather}
where the Hamiltonian \mb{\cR_{D}} is
\begin{gather}
\cR_{D} := \cR+\nabla(r\!-\!\omegaone)
+\frac{1}{2i\hbar}[r\!-\!\omegaone\ocomma r\!-\!\omegaone]
 = \cR+(\nabla-\delta)r+\frac{1}{2i\hbar}[r\ocomma r]-\omega .
\label{hamcrdcurv}
\end{gather}
In the third and fourth equality of \eq{rdcurv} we use that \mb{\nabla}
respects the \mb{\circ} product \e{nablarespectcirc} and the Jacobi identity
for the \mb{\circ} product, respectively. In the second equality of
\eq{hamcrdcurv} we use that \mb{(\nabla\omegaone)=0},
\mb{\delta=\{\omegaone,\cdot\}} and \mb{\{\omegaone,\omegaone\}=-2\omega}.

\section[Flat/nilpotent $D$ connection]{Flat/nilpotent $\boldsymbol{D}$ connection}
\label{secflatdconn}

The next main principle of Fedosov quantization is to choose the \mb{D}
connection to be f\/lat, or equivalently, nilpotent:
\begin{gather}
R_{D} \equiv D^{2} \equiv \Hf[D,D] = 0 .\label{rzerocurvature}
\end{gather}
In other words, the odd \mb{D} operator is a deformation of the odd Koszul--Tate
dif\/ferential \mb{\delta}, such that the nilpotency is preserved, cf.\
\eqs{deltanilp}{rzerocurvature}. (This setup is similar to the construction of
an odd nilpotent BRST operator
\cite{fravil75,batvil77,frafra78,batfra83a,batfra83b,fhst89,hentei92}.) Since
we want to achieve the nilpotency \e{rzerocurvature}, it now becomes clear why
it was so important that the \mb{x}- and \mb{c}-derivative terms in the \mb{R}
curvature two-form \e{curvasham} vanish. This is because there are no other
\mb{x}- and \mb{c}-derivatives in the \mb{R_{D}} curvature two-form \e{rdcurv}
to cancel them. (All other derivatives in \eq{rdcurv} are \mb{y}-derivatives.)
This crucial point is the main reason that the \mb{\nabla} connection is
assumed to be torsion-free. (References~\cite{karsch01,kar03} consider a hybrid
model where torsion is allowed in the \mb{y}-sector but not in the
\mb{c}-sector to avoid the \mb{x}- and \mb{c}-derivatives. Note, however, that
they restrict the possible torsion by imposing both the independent conditions
\mb{(\nabla_{i}\tilde{\omega}_{jk})=0} and \mb{(d\omega)=0} at the same time.)
{}For the curvature two-form \mb{R_{D}} to be zero, it is enough to let the
Hamiltonian curvature two-form \mb{\cR_{D}} be {\em Abelian}, i.e.,  to let it
belong to the center
\begin{gather*}
Z(\cA) := \{a\in\cA \mid [\cA \ocomma a]=0\}
 = \{a\in\cA \mid \{\cA,a\}=0\} = \cA_{\bullet 0}
 \equiv \Gamma\left(\bigwedge{}^{\bullet}(T^{*}M)\right)[[\hbar]] 
\end{gather*}
of the algebra \mb{(\cA,+,\circ)}. In other words, there should exist an
Abelian two-form \mb{\casimir\in\cA_{20}}, such that
\begin{gather}
 \cR_{D} = -\casimir-\omega \in \cA_{20} .\label{rabelian}
\end{gather}
Here \mb{\omega} is just the symplectic two-form \e{omegatwoform} itself.
(Recall that \mb{\omega} trivially belongs to the center~\mb{Z(\cA)}. The
signs and the shift with \mb{\omega} in \eq{rabelian} are introduced without
loss of generality for later convenience.) The Abelian condition \e{rabelian}
turns into Fedosov's \mb{r}-equation
\begin{gather}
(\delta r) = \cR+\casimir+(\nabla r) +\frac{1}{2i\hbar}[r\ocomma r] .
\label{fedosovreq}
\end{gather}

\begin{theorem}[Fedosov's 1st theorem]\label{fedosov1theorem}
Let there be given an Abelian two-form \mb{\casimir\in\cA_{20}} that is
closed $($\mb{\equiv}symplectic$)$,
\begin{gather*}
(d\casimir) = 0 , 
\end{gather*}
and that satisfies the boundary condition
\begin{gather}
\casimir_{(0)} = 0 .
\label{bigomegabc}
\end{gather}
Then there exists a unique one-form valued \mb{r}-solution
\begin{gather*}
r = \oplus_{\nn\geq 0}r_{(\nn)} \in \cA_{1\bullet} ,
\qquad    \deg(r_{(\nn)}) = \nn ,
\end{gather*}
to Fedosov's equation \eqref{fedosovreq} such that \mb{r} is \mb{\delta^{*}}-closed,
\begin{gather}
(\delta^{*}r) = 0 ,
\label{fedosov1star}
\end{gather}
and satisfies the boundary condition
\begin{gather}
r_{(1)} = 0 \label{r1zero}
\end{gather}
for the \mb{(1)}-sector. As a consequence, it turns out that the first three
sectors \mb{r_{(0)}}, \mb{r_{(1)}} and \mb{r_{(2)}} are identically zero.
\end{theorem}

\begin{proof}
Let us split the Abelian condition \e{fedosovreq} in Fedosov degree:
\begin{gather}
\delta r_{(0)} \equiv 0 ,\label{r0} \\
\delta r_{(\nn+1)} = \cR_{(\nn)}+\casimir_{(\nn)}
+\sum_{k=0}^{\nn}\stackrel{(\nn-k)}{\nabla} r_{(k)} \nonumber\\
\phantom{\delta r_{(\nn+1)} =}{} +\frac{1}{2i\hbar}\sum_{0\leq k,\ell}^{k+\ell\leq \nn+2}
[r_{(k)}\ocomma r_{(\ell)}]_{(\nn+2-k-\ell)} \qquad {\rm for} \quad \nn\geq 0
 . \label{rm1}
\end{gather}
In \eq{rm1} the Hamiltonian curvature two-form \mb{\cR} and the closed Abelian
two-form \mb{\casimir} have also been expanded in Fedosov degree
\begin{alignat*}{4}
& \cR=\oplus_{\nn\geq 2}\cR_{(\nn)} \in \cA_{22},\qquad &&
\deg(\cR_{(\nn)})=\nn,\qquad && \cR_{(\odd)}=0 , & \nonumber\\
& \casimir= \oplus_{\nn\geq 0}\casimir_{(\nn)} \in \cA_{20} , \qquad &&
\deg(\casimir_{(\nn)}) = \nn , \qquad && \casimir_{(\odd)} = 0 . & 
\end{alignat*}
A priori it is known that the \mb{(0)}-sector \mb{r_{(0)}=c^{i}\eta_{i}(x)} is
a \mb{y}- and \mb{\hbar}-independent one-form. The \eq{r0} is therefore
automatically satisf\/ied. It follows from \eq{fedosov1star} that
\begin{gather*}
0 = \delta^{*}r_{(0)} = y^{i}\eta_{i}(x) .
\end{gather*}
Therefore the \mb{(0)}-sector
\begin{gather*}
r_{(0)} = 0 
\end{gather*}
vanishes identically. Equation \e{rm1} with \mb{\nn=0} is automatically
satisf\/ied  because of the two boundary conditions \es{bigomegabc}{r1zero}.
Putting \mb{\nn=1} in \eq{rm1} yields \mb{\delta r_{(2)}=0}. Hence the
\mb{(2)}-sector \mb{r_{(2)}} is a one-form that is both \mb{\delta}-closed and
\mb{\delta^{*}}-closed, and therefore it must be identically zero as well:
\begin{gather*}
r_{(2)} = 0 .
\end{gather*}
Since \mb{r_{(0)}}, \mb{r_{(1)}} and \mb{r_{(2)}} are zero, the right-hand side  expression
for \mb{\delta r_{(\nn+1)}} in \eq{rm1} will only depend on previous entries
\mb{r_{(\leq \nn)}}. Hence \eq{rm1} is a recursive relation. The consistency
relation for the Abelian condition \e{fedosovreq} is that the right-hand side  should be
\mb{\delta}-closed. This is indeed the case:
\begin{gather}
\delta(\RHS) = \delta\left(\cR+\casimir+\nabla r
+\frac{1}{2i\hbar}[r\ocomma r]\right)
 = \delta\cR+[\delta,\nabla]r-\nabla(\delta r)
+\frac{1}{i\hbar}[\delta r\ocomma r] = -D(\delta r) \nonumber\\
\phantom{\delta(\RHS)}{} =-\nabla\!\left(\cR+\casimir+\nabla r +\frac{1}{2i\hbar}[r\ocomma r]\right)\!
 +\frac{1}{i\hbar}\left[\cR+\casimir+\nabla r +\frac{1}{2i\hbar}[r\ocomma r]
\ocomma r\right] = 0 .\!\!\!\!\!
\label{crr}
\end{gather}
In the second equality of \eq{crr} we have used that \mb{\delta} respects the
\mb{\circ} product. In the third equality we have used that \mb{(\delta\cR)=0},
that \mb{[\delta,\nabla]=0}, and that \mb{\delta} is nilpotent. In the f\/ifth
(=last) equality we have used that \mb{(\nabla\cR)=0}, that \mb{(d\casimir)=0},
that \mb{\nabla^{2}=\{\cR,\cdot\}}, that \mb{\nabla} respects the \mb{\circ}
product, and the Jacobi identity for the \mb{\circ} product.

We now prove by induction on the Fedosov degree \mb{(\nn)} that there exists a
unique solution~\mb{r_{(\nn+1)}} to \eq{rm1} if there exists a unique solution
for all the previous entries~\mb{r_{(\leq \nn)}}. This is essentially a
consequence of the Poincar\'e lemma~\ref{poincarelemma}. Uniqueness:  The
dif\/ference
\begin{gather*}
\Delta r_{(\nn+1)} := r_{(\nn+1)}^{\prime} - r_{(\nn+1)}^{\prime\prime}
\end{gather*}
between two solutions \mb{r_{(\nn+1)}^{\prime}} and
\mb{r_{(\nn+1)}^{\prime\prime}} must satisfy the homogeneous version
\mb{\delta(\Delta r_{(\nn+1)})=0} of \eq{rm1}, i.e.,  with no right-hand side.
Hence the dif\/ference \mb{\Delta r_{(\nn+1)}} is a one-form that is both
\mb{\delta}-closed and \mb{\delta^{*}}-closed, and therefore it must be
identically zero. Existence:  Def\/ine
\begin{gather*}
r_{(\nn+1)} := \delta^{-1}(\RHS_{(\nn)}) \qquad {\rm for} \quad \nn\geq 0 ,
\end{gather*}
where \mb{\RHS_{(\nn)}} is the two-form valued right-hand side  of \eq{rm1}. This clearly
def\/ines a \mb{\delta^{*}}-closed one-form \mb{r_{(\nn+1)}}. To check \eq{rm1},
it is enough to check that the two-form \mb{\RHS_{(\nn)}} is
\mb{\delta}-closed. But this follows by linearity from the consistency relation
\e{crr}, because~\mb{\nabla} and \mb{\circ} both carry positive Fedosov
degree, and the f\/irst three \mb{r}-sectors vanish to cancel the negative
{}Fedosov degree coming from the \mb{\hbar^{-1}}-factor, so that only previous
entries \mb{r_{(\leq \nn)}} can participate to the \mb{(\nn)}-sector.
\end{proof}

We emphasize that the unique deformation one-form \mb{r} from Fedosov's 1st
theorem is glo\-bally well-def\/ined, since it basically appeared from inverting
the Koszul--Tate \mb{\delta} dif\/ferential, cf.\  Section~\ref{secdeltam1}.
Normally, one would choose a trivial Abelian two-form \mb{\casimir\equiv 0}.
Also note that the two-form \mb{\cR + \casimir} is the lone source term that
forces \mb{r} to be non-trivial. We list here the f\/irst few unique
\mb{r}-terms:
\begin{gather*}
r_{(0)} = 0 ,  \qquad r_{(1)} = 0 ,  \qquad r_{(2)} = 0 , \qquad
r_{(3)} = \delta^{-1}\left(\cR_{(2)}+\casimir_{(2)}\right) ,\nonumber\\
r_{(4)} = \delta^{-1}\big(\stackrel{(0)}{\nabla} r_{(3)}\big) ,  \qquad \ldots.
\end{gather*}
Similarly, the f\/irst few terms in the \mb{D} expansion read
\begin{gather*}
\stackrel{(-1)}{D} = -\delta ,\qquad
\stackrel{(0)}{D} = \stackrel{(0)}{\nabla} ,\qquad
\stackrel{(1)}{D} = \frac{1}{i\hbar}[r_{(3)}\ocomma\cdot\, ]_{(0)} , \qquad  \ldots.
\end{gather*}

\section{Horizontal sections}
\label{sechorizontalsec}

Fedosov's 1st theorem establishes the existence of a globally well-def\/ined,
unique, f\/lat/nilpo\-tent~\mb{D} connection. Since this higher-order connection~\mb{D} is f\/lat, it is possible to solve the horizontal condition \mb{(Da)=0}
locally for a zero-form valued section \mb{a\in\cA_{0\bullet}}. (In other
words, the f\/latness relation \e{rzerocurvature} is the local consistency
relation for the horizontal condition.) As we shall see below there is no
obstruction in patching together local horizontal sections \mb{a} into global
horizontal sections, basically because \mb{D} is a deformation of (minus) the
Koszul--Tate \mb{\delta}-dif\/ferential.

\begin{theorem}[Fedosov's 2nd theorem]\label{fedosov2theorem}
Let there be given a quantum function $($also known as  a symbol$)$
\mb{f\in C^{\infty}(M)[[\hbar]]\equiv\cA^{}_{00}}.
Then there exists a unique zero-form valued section
\begin{gather*}
a = \oplus_{\nn\geq 0}a_{(\nn)} \in \cA_{0\bullet}
 ,\qquad \deg(a_{(\nn)}) = \nn ,
\end{gather*}
that is \mb{D}-horizontal
\begin{gather*}
(Da) = 0 ,
\end{gather*}
and that satisfies the boundary condition
\begin{gather}
a|^{}_{y=0} \equiv a^{}_{00} = f .
\label{fedosov2cond}
\end{gather}
\end{theorem}

\begin{proof}
First note that that a
zero-form \mb{a\in\cA_{0\bullet}} is automatically \mb{\delta^{*}}-closed:
\begin{gather*}
(\delta^{*}a) \equiv 0 ,\qquad  a\in\cA_{0\bullet} .
\end{gather*}
The horizontal condition \mb{(Da)=0} becomes
\begin{gather}
(\delta a) = (\nabla a)+\frac{1}{i\hbar}[r\ocomma a] .
\label{fedosovaeq}
\end{gather}
Let us split the horizontal condition \e{fedosovaeq} in Fedosov degree:
\begin{gather}
\delta a_{(0)} \equiv 0 ,\nonumber \\ 
\delta a_{(\nn+1)} = \sum_{k=0}^{\nn}\stackrel{(\nn-k)}{\nabla} a_{(k)}
+\frac{1}{i\hbar}\sum_{0\leq k,\ell}^{k+\ell\leq \nn+2}
[r_{(k)}\ocomma a_{(\ell)}]_{(\nn+2-k-\ell)}\qquad {\rm for} \quad \nn\geq 0 .
\label{am1}
\end{gather}
Note that the right-hand side  expression for \mb{\delta a_{(\nn+1)}} only depends on
previous entries~\mb{a_{(\leq \nn)}}, because \mb{r_{(0)}}, \mb{r_{(1)}} and
\mb{r_{(2)}} are zero. Hence \eq{am1} is a recursive relation.
The consistency relation for the horizontal condition \e{fedosovaeq} is that
the right-hand side  should be~\mb{\delta}-closed. This is indeed the case:
\begin{gather}
\delta(\RHS) = \delta\left(\nabla a+\frac{1}{i\hbar}[r\ocomma a]\right)
 = [\delta,\nabla]a-\nabla(\delta a)-\frac{1}{i\hbar}[ r\ocomma \delta a]
+\frac{1}{i\hbar}[\delta r\ocomma a] \nonumber\\
 \phantom{\delta(\RHS)}{} = -\nabla\left(\nabla a+\frac{1}{i\hbar}[r\ocomma a]\right)
-\frac{1}{i\hbar}\left[r\ocomma\nabla a
+\frac{1}{i\hbar}[r\ocomma a] \right] \nonumber\\
 \phantom{\delta(\RHS)=}{} +\frac{1}{i\hbar}\left[\cR+\casimir+\nabla r +\frac{1}{2i\hbar}[r\ocomma r]
\ocomma a \right] = 0 .\label{cra}
\end{gather}
In the second equality of \eq{cra} we have used that \mb{\delta} respects the
\mb{\circ} product. In the third equality we have used that
\mb{[\delta,\nabla]=0}. In the fourth (=last) equality we have used that
\mb{\nabla^{2}=\{\cR,\cdot\}}, that \mb{\nabla} respects the \mb{\circ}
product, and the Jacobi identity for the \mb{\circ} product.

We now prove by induction on the Fedosov degree \mb{(\nn)} that there exists a
unique solution \mb{a_{(\nn+1)}} to \eq{am1} if there exists a unique solution
for all the previous entries \mb{a_{(\leq \nn)}}. This is essentially a
consequence of the Poincar\'e lemma~\ref{poincarelemma}.
Uniqueness:  The dif\/ference
\begin{gather*}
\Delta a_{(\nn+1)} := a_{(\nn+1)}^{\prime} - a_{(\nn+1)}^{\prime\prime}
\end{gather*}
between two solutions \mb{a_{(\nn+1)}^{\prime}} and
\mb{a_{(\nn+1)}^{\prime\prime}} must satisfy the homogeneous versions
\mb{\delta(\Delta a^{}_{(\nn+1)})=0} and \mb{\Delta a^{}_{(\nn+1)00}=0} of the
horizontal condition \e{am1} and boundary condition \e{fedosov2cond}, i.e.,  with
no right-hand sides. Hence the dif\/ference \mb{\Delta a_{(\nn+1)}} is both
\mb{\delta}-closed, \mb{\delta^{*}}-closed and with no \mb{00}-sector.
Therefore it must be identically zero. Existence:  Def\/ine
\begin{gather*}
a_{(0)} := f_{(0)} , \\ 
a_{(\nn+1)} := f_{(\nn+1)}+\delta^{-1}(\RHS_{(\nn)})\qquad {\rm for} \quad \nn\geq 0 ,
\end{gather*}
where \mb{\RHS_{(\nn)}} is the one-form valued right-hand side  of \eq{am1}. This clearly
def\/ines a~zero-form \mb{a_{(\nn+1)}} that satisf\/ies the boundary condition
\e{fedosov2cond}. To check \eq{am1}, it is enough to check that the one-form
\mb{\RHS_{(\nn)}} is~\mb{\delta}-closed. But this follows by linearity from the
consistency relation~\e{cra}, because~\mb{\nabla} and~\mb{\circ} both carry
positive Fedosov degree, and the f\/irst three~\mb{r}-sectors vanish to cancel
the negative Fedosov degree coming from the \mb{\hbar^{-1}}-factor, so that
only previous entries \mb{a_{(\leq \nn)}} can participate to the
\mb{(\nn)}-sector.
\end{proof}

We list here the solution to the unique f\/irst-order correction \mb{a_{(1)}}:
\begin{gather*}
a_{(1)}-f_{(1)} = \delta^{-1}(\stackrel{(0)}{\nabla} a_{(0)})
 = y^{i}\bigg(\papal{x^{i}}f_{(0)}\bigg) ,
\end{gather*}
which we'll use in the next Section~\ref{secstarproduct}.

\section[The $*$ product]{The $\boldsymbol{*}$ product}
\label{secstarproduct}

Fedosov's 2nd theorem establishes an isomorphism
\begin{gather}
(\cA^{}_{00},+,*) \ni f \ \stackrel{Q}{\mapsto}  \ Q(f) \in (W^{}_{D},+,\circ)
\label{qmap}
\end{gather}
between the algebra \mb{\cA^{}_{00}} of quantum functions (=symbols),
\begin{gather*}
\cA^{}_{00} = C^{\infty}(M)[[\hbar]]
 = \{a\in\cA \mid (\delta a)=0=(\delta^{*}a)\} ,
\end{gather*}
and the algebra of zero-form valued horizontal sections,
\begin{gather*}
W^{}_{D} := \{ a\in\cA^{}_{0\bullet} \mid (Da)=0 \} . 
\end{gather*}
That the vector space \mb{(W^{}_{D},+,\circ)} is an subalgebra, i.e.,  closed
with respect to the  \mb{\circ} product, follows basically because the connection \mb{\nabla}
respects the \mb{\circ} product, cf.\  \eq{nablarespectcirc}. The \mb{*} product
in diagram \e{qmap} is by def\/inition induced from the \mb{\circ} product as
\begin{gather*}
f*g := Q^{-1}\left(Q(f)\circ Q(g)\right) .
\end{gather*}
This \mb{*} product obviously inherits associativity from the \mb{\circ}
product, and the \mb{Q} map \e{qmap} is obviously an algebra isomorphism.
Moreover, the inverse map
\begin{gather*}
(W^{}_{D},+,\circ) \ni a \stackrel{Q^{-1}}{\mapsto}
a|^{}_{y=0} \equiv a^{}_{00} \in (\cA^{}_{00},+,*)
\end{gather*}
is simply the restriction to \mb{y\!=\!0}, cf.\  boundary condition
\e{fedosov2cond}. It remains to check that the classical boundary condition
\e{classicalbc0} holds. Let us f\/irst expand to the appropriate orders:
\begin{gather*}
Q(f) = a_{(0)}+a_{(1)}+a_{(2)}+{\cal O}(\hbar^{2},y\hbar,y^{3}) , \\
Q(g) = b_{(0)}+b_{(1)}+b_{(2)}+{\cal O}(\hbar^{2},y\hbar,y^{3}) , \\
Q(f)\circ Q(g) = Q(f) Q(g)
+\frac{i\hbar}{2}Q(f)\papar{y^{j}}\mult^{jk}\papal{y^{k}}Q(g)
+{\cal O}(\hbar^{2}) .
\end{gather*}
Therefore
\begin{gather}
\lim_{\hbar\to 0} f*g = f_{(0)}g_{(0)} , \label{classicalbc00} \\
\lim_{\hbar\to 0} \frac{f*g-f_{(0)}g_{(0)}}{i\hbar}
 = \frac{a_{(2)}|^{}_{y=0}}{i\hbar}b_{(0)}
+a_{(0)}\frac{b_{(2)}|^{}_{y=0}}{i\hbar}
+\Hf a_{(1)}\papar{y^{j}}\mult^{jk}_{(0)}\papal{y^{k}}b_{(1)} ,\nonumber
\end{gather}
or equivalently,
\begin{gather}
\lim_{\hbar\to 0} \frac{f*g-fg}{i\hbar}
 = \Hf f_{(0)}\papar{x^{j}}\mult^{jk}_{(0)}\papal{x^{k}}g_{(0)} .
\label{classicalbc01}
\end{gather}
Equations \es{classicalbc00}{classicalbc01} are precisely the content of the
classical boundary condition \e{classicalbc0}. It is also easy to check
condition \e{starunit}.

\medskip

\noindent
{\bfseries\itshape Note added.}
The author has kindly been informed by a referee of the existence of the
Ph.D.\ Thesis \cite{eckel00}, which the author has been unable to obtain, and which
considers Fedosov quantization on supermanifolds.

\subsection*{Acknowledgement}

The author thanks I.A.~Batalin, D.~Sternheimer and the three referees for
comments. The work of K.B.\ is supported by the Ministry of Education of the
Czech Republic under the project MSM 0021622409.

\bigskip

\noindent
{\bfseries\itshape Editorial Comments.} This paper presents the Fedosov construction in a fairly general
framework. It is generally known that the tensor need not be skew-symmetric
as long as its skew-symmetric part is non-degenerate, can be a series in
the deformation parameter, and that Fedosov's construction has also been
extended to the case of super-manifolds. However it appears from the
last reports that there is a point in publishing the unif\/ied
presentation of the author, which should be viewed as a good review of
Fedosov's construction in a quite general context, even if the author
overemphasizes somewhat the importance of the generalizations he considers.

\pdfbookmark[1]{References}{ref}
\LastPageEnding

\end{document}